\documentclass{article}
\usepackage{graphicx}
\usepackage{amsmath}

\newtheorem{theorem}{Theorem}

\newenvironment{proof}[1][Proof]{\textbf{#1.} }{\ \rule{0.5em}{0.5em}}
\input{tcilatex}

\begin{document}

\title{Computation of the Generalized $F$ Distribution}
\author{\begin{tabular}{ll}
Charles F. Dunkl & Donald E. Ramirez
\end{tabular}
\\
%EndAName
Department of Mathematics\\
University of Virginia\\
Charlottesville, VA 22903-3199 USA}
\date{}
\maketitle

\begin{abstract}
Exact expressions for the distribution function of a random variable of the
form $((\alpha _{1}\chi _{m_{1}}^{2}+\alpha _{2}\chi _{m2}^{2})/|m|)/(\chi
_{\nu }^{2}/\nu )$ are given where the chi-square distributions are
independent with degrees of freedom $m_{1},m_{2},$ and $\nu $ respectively.
Applications to detecting joint outliers and Hotelling's misspecified $T^{2}$
distribution are given.
\end{abstract}


\section{\textbf{Introduction}}

The \textit{generalized }$F$\textit{\ distribution }is defined as follows.
Suppose that the elements of $\mathbf{X}=[\chi _{m_{1}}^{2},\cdots ,\chi
_{m_{r}}^{2}]^{^{\prime }}$ $(r>1)$ are independent chi-square random
variables with degrees of freedom $(m_{1},\cdots ,m_{r}),$ respectively; let 
$\{\alpha _{1}\geq \alpha _{2}\geq \cdots \geq \alpha _{r}>0\}$ be
nonincreasing positive weights; and identify $T=\alpha _{1}\chi
_{m_{1}}^{2}+\cdots +\alpha _{r}\chi _{m_{r}}^{2}.$ If $\mathcal{L}(V)=\chi
^{2}(\nu )$ independently of $\mathbf{X}$, then the $cdf$ of 
\begin{equation}
W=\frac{T/|m|}{V/\nu }=\frac{(\alpha _{1}\chi _{m_{1}}^{2}+\cdots +\alpha
_{r}\chi _{m_{r}}^{2})/|m|}{V/\nu },
\end{equation}
where $|m|=m_{1}+\cdots +m_{r},$ is denoted by $F_{r}(w;\alpha _{1},\cdots
,\alpha _{r};m_{1},\cdots ,m_{r};\nu )$. If all of the $\alpha _{i}$ $\left(
1\leq i\leq r\right) $ are equal to say $\alpha ,$ then the $cdf$ of $W$ is
denoted by $F_{r}(w;\alpha ;m_{1},\cdots ,m_{r};\nu )$, the scaled central $%
F $ distribution with degrees of freedom $(|m|,\nu )$. To avoid the trivial
case, we will assume that the positive weights are pairwise distinct.

We will give exact expressions for the $pdf$ of $W$ for $r=2$ in terms of
the hypergeometric series $_{2}F_{1}.$ This is the analog for generalized
functions of the known result for a mixture of two chi-square distributions
(Bock and Solomon (1988))$.$ For $r>2,$ we give three numerically tractable
expressions for the $pdf$ and $cdf$ of $W$. Applications include the
detection of joint outliers using Cook's $D_{I}$ statistics and the
calculation of the power of Hotelling's $T^{2}$ test with a misspecifed
scale.

\section{\textbf{The Distribution of }$(T/|m|)/(V/\protect\nu )$}

Building on the work of Robbins and Pitman (1949), Gurland (1955), and Kotz,
Johnson, and Boyd (1967), Ramirez and Jensen (1991) showed how to compute
the $pdf$ for $W_{0}=T/V$ as a weighted series of $F$ distributions; and
they computed the error bounds for the truncated partial sums. Their results
are stated for $W_{0}=T/V,$ with $r=p,$ and with $\mathcal{L}(V)=\chi
^{2}(\nu -p+1);$ and they used the notation from Kotz, Johnson and Boyd
(1967). We give the results for the general case below where it is
convenient for our derivation to use the notation from Robbins and Pitman
(1949).

\subsection{The Probability Distribution Function for $W$}

Write 
\begin{equation}
T=\alpha _{r}(\frac{\alpha _{1}}{\alpha _{r}}\chi _{m_{1}}^{2}+\cdots +\frac{%
\alpha _{r-1}}{\alpha _{r}}\chi _{m_{r-1}}^{2}+\chi _{m_{r}}^{2}).
\end{equation}
Following Robbins and Pitman (1949, p. 555) define the constants $c_{j}$ by
the identity 
\begin{equation}
A\,\prod_{i=1}^{r}\left( 1-u_{i}z\right) ^{-m_{i}/2}=\sum_{j=0}^{\infty
}c_{j}z^{j},  \label{R&P}
\end{equation}

where 
\begin{equation}
A=\prod_{i=1}^{r}\left( \frac{\alpha _{i}}{\alpha _{r}}\right) ^{-m_{i}/2}%
\text{.}
\end{equation}
The series in Equation \ref{R&P} converges absolutely for $|z|<\alpha
_{1}/(\alpha _{1}-\alpha _{r})$. Set $z=0$ to see that $c_{0}=A$, and set $%
z=1$ for the equality $\sum_{j=0}^{\infty }c_{j}=1$. Then $P[T\leq y]=\sum
c_{j}G_{|m|+2j}(y/\alpha _{r})$, where $G_{k}$ is the $cdf$ for the
chi-square distribution with $k$ degrees of freedom. As in Ramirez and
Jensen (1991, p. 100), we find that the $pdf$ for $W=(T/|m|)/(V/\nu )$ has
the representation as stated in the following

\begin{theorem}
With the notation above, 
\begin{eqnarray}
h_{W}(w) &=&\sum_{j=0}^{\infty }\frac{|m|}{\nu }\frac{c_{j}}{\alpha _{r}}%
\frac{\Gamma \left( \frac{\nu +|m|+2j}{2}\right) \left( \frac{|m|}{\nu }%
\frac{w}{\alpha _{r}}\right) ^{\left( |m|+2j-2\right) /2}}{\Gamma \left( 
\frac{|m|+2j}{2}\right) \Gamma \left( \frac{\nu }{2}\right) \left( 1+\frac{%
|m|}{\nu }\frac{w}{\alpha _{r}}\right) ^{\left( \nu +|m|+2j\right) /2}} 
\notag \\
&=&\sum_{j=0}^{\infty }\frac{c_{j}}{\alpha _{r}}\frac{|m|}{|m|+2j}%
\,\,f_{F}\left( \frac{|m|}{|m|+2j}\frac{w}{\alpha _{r}};|m|+2j,\nu \right) ,
\label{W}
\end{eqnarray}
with $f_{F}(w;v_{1},v_{2})$ the density of the central $F$ distribution with
degrees of freedom $(v_{1},v_{2}).$ \newline
A bound for the global truncation error $e_{\tau }$ for the $\tau ^{th}$
partial sum of the $pdf$ of $W=(T/|m|)/(V/\nu )$ is given by 
\begin{eqnarray}
&&\sum_{j=\tau +1}^{\infty }\frac{c_{j}}{\alpha _{r}}\frac{|m|}{|m|+2j}%
\,\,f_{F}\left( \frac{|m|}{|m|+2j}\frac{w}{\alpha _{r}};|m|+2j,\nu \right) \\
&\leq &\frac{|m|}{\alpha _{r}(|m|+2(\tau +1))}(1-(c_{0}+\cdots +c_{\tau
}))=e_{\tau }.  \label{gWerr}
\end{eqnarray}
\end{theorem}

\begin{proof}
Use the equality $\sum_{i=0}^{\infty }c_{i}=1,$ and note that $\left|
f_{F}(w;v_{1},v_{2})\right| \leq 1$ when $v_{1}\geq 2$ and $v_{2}\geq 1$.
\end{proof}

The global bound $e_{\tau }$ can be used to determine the number of terms $%
\tau $ to use in the truncated series expansion of the $pdf$ for $W$ in
Equation \ref{W}. In Section \ref{local bound pdf}, we improve on the global
error bound $e_{\tau }$ by identifying the local error bound as a
hypergeometric function $_{2}F_{1}.$

\subsection{Calculation of the Coefficients $c_{j}$}

Kotz, Johnson, and Boyd (1967) gave the following expression for $c_{j},$

\begin{eqnarray}
c_{0} &=&\prod_{i=1}^{r}\left( \frac{\alpha _{r}}{\alpha _{i}}\right)
^{m_{i}/2}=A,  \notag \\
d_{j} &=&\sum_{i=1}^{r}\frac{m_{i}}{2}\left( 1-\frac{\alpha _{r}}{\alpha _{i}%
}\right) ^{j},\,\,j\geq 1, \\
c_{j} &=&\frac{1}{j}\sum_{l=0}^{j-1}\left( d_{j-l}c_{l}\right) ,\,\,j\geq 1.
\notag
\end{eqnarray}

We are able to reduce the numerical complexity in the computation of the
coefficients $c_{j}$ by determining a recursive algorithm for $c_{j}.$ Fix
parameters $\mu _{1},\ldots ,\mu _{r}$ and variables $u_{1},\ldots ,u_{r}$
with $|u_{i}|<1$ for all $i\,(1\leq i\leq r).$ For $k=0,1,2,..$., let 
\begin{equation*}
P_{k}=\sum_{|\mathbf{n}|=k}\prod_{i=1}^{r}\frac{(\mu _{i})_{n_{i}}}{n_{i}!}%
u_{i}^{n_{i}},\mathbf{n}=(n_{1},\ldots ,n_{r}).
\end{equation*}
Note that $\sum_{k=0}^{\infty }P_{k}=\prod_{i=1}^{r}(1-u_{i})^{-\mu _{i}}$.
Denote the set $R=\{1,2,\ldots ,r\}.$ For $i\in R,$ define 
\begin{eqnarray*}
e_{i} &=&\sum_{S\subset R,|S|=i}\prod_{j\in S}u_{j}, \\
f_{i} &=&\sum_{S\subset R,|S|=i}\left( \sum_{j\in S}\mu _{j}\right)
\prod_{j\in S}u_{j}.
\end{eqnarray*}
Thus $e_{i}$ is the elementary symmetric function of degree $i$ in $%
u_{1},\ldots ,u_{r}.$ Then for $k\geq 1$%
\begin{equation*}
kP_{k}=\sum_{i=1}^{r}(-1)^{i-1}((k-i)e_{i}+f_{i})P_{k-i}.
\end{equation*}
To prove the identity, let $\lambda _{i}=\mu _{i}-1$ for all $i;$ and for a
fixed $\mathbf{n}=(n_{1},\ldots ,n_{r})$ with $|\mathbf{n}|=k,$ examine the
coefficient of $\dprod_{i=1}^{r}\dfrac{(\mu _{i})_{n_{i}}}{n_{i}!}%
u_{i}^{n_{i}}$ \ in the sum $kP_{k}+%
\sum_{i=1}^{r}(-1)^{i}((k-i)e_{i}+f_{i})P_{k-i}.$ Let $\zeta _{i}=\dfrac{%
n_{i}}{\mu _{i}+n_{i}-1}=\dfrac{n_{i}}{\lambda _{i}+n_{i}}$ then this
coefficient equals 
\begin{equation*}
k+\sum_{i=1}^{r}(-1)^{i}\sum_{S\subset R,|S|=i}(k+\sum_{j\in S}\lambda
_{j})\prod_{j\in S}\zeta _{j}.
\end{equation*}
The coefficient of $k$ in this expression is $\prod_{i=1}^{r}(1-\zeta _{i}).$
For each $s,$ the coefficient of $\lambda _{s}$ is 
\begin{gather*}
\zeta _{s}\sum_{i=1}^{r}(-1)^{i}\sum \left\{ \prod_{j\in S}\zeta
_{j}:S\subset R\backslash \{s\},|S|=i-1\right\} \\
=-\zeta _{s}\prod_{i\neq s}(1-\zeta _{i}).
\end{gather*}
But $\lambda _{s}\zeta _{s}=\dfrac{\lambda _{s}n_{s}}{\lambda _{s}+n_{s}}%
=n_{s}(1-\zeta _{s}),$ and so these terms sum to $-\sum_{s=1}^{r}n_{s}%
\prod_{i=1}^{r}(1-\zeta _{i}),$ and $|\mathbf{n}|=k.$ This completes the
proof by noting that $c_{k}=AP_{k}$ with $\mu _{i}=m_{i}/2$.

\section{Exact Expressions for the $pdf$ of $W$}

Use the negative binomial series 
\begin{equation}
(1-sz)^{-b}=\sum_{m=0}^{\infty }s^{m}\frac{(b)_{m}}{m!}z^{m}
\end{equation}
to express Equation \ref{R&P} as 
\begin{equation}
\sum_{j=0}^{\infty }c_{j}z^{j}=A\sum_{j=0}^{\infty }z^{j}\sum_{i_{1}+\cdots
+i_{r-1}=j}\prod_{k=1}^{r-1}\frac{u_{k}^{i_{k}}}{i_{k}!}\left( \frac{m_{k}}{2%
}\right) _{i_{k}},
\end{equation}
with 
\begin{equation*}
0\leq u_{i}=1-\frac{\alpha _{r}}{\alpha _{i}}<1\,(1\leq i\leq r).
\end{equation*}
Note that $u_{r}=0$. Denote 
\begin{eqnarray}
a &=&\frac{\nu \alpha _{r}}{|m|}, \\
B_{0} &=&a^{\nu /2}\frac{\Gamma \left( \frac{\nu +|m|}{2}\right) }{\Gamma
\left( \frac{|m|}{2}\right) \Gamma \left( \frac{\nu }{2}\right) }, \\
B_{1}(w) &=&\frac{w^{\left( |m|-2\right) /2}}{\left( a+w\right) ^{\left( \nu
+|m|\right) /2}},
\end{eqnarray}
and write the $pdf$ for $W=(T/|m|)/(V/\nu )$ with 
\begin{equation*}
t(w)=\frac{w}{a+w},
\end{equation*}
as 
\begin{gather}
h_{W}(w)=\sum_{j=0}^{\infty }\frac{|m|}{\nu }\frac{c_{j}}{\alpha _{r}}\frac{%
\Gamma \left( \frac{\nu +|m|+2j}{2}\right) \left( \frac{|m|}{\nu }\frac{w}{%
\alpha _{r}}\right) ^{\left( |m|+2j-2\right) /2}}{\Gamma \left( \frac{|m|+2j%
}{2}\right) \Gamma \left( \frac{\nu }{2}\right) \left( 1+\frac{|m|}{\nu }%
\frac{w}{\alpha _{r}}\right) ^{\left( \nu +|m|+2j\right) /2}}  \notag \\
=B_{0}B_{1}(w)\sum_{j=0}^{\infty }\frac{c_{j}\left( \frac{\nu +|m|}{2}%
\right) _{j}}{\left( \frac{|m|}{2}\right) _{j}}t(w)^{j}  \label{hw line 2} \\
=AB_{0}B_{1}(w)\sum_{j=0}^{\infty }\frac{\left( \frac{\nu +|m|}{2}\right)
_{j}}{\left( \frac{|m|}{2}\right) _{j}}t(w)^{j}\sum_{i_{1}+\cdots
+i_{r-1}=j}\prod_{k=1}^{r-1}\frac{u_{k}^{i_{k}}}{i_{k}!}\left( \frac{m_{k}}{2%
}\right) _{i_{k}}  \label{Wline3} \\
=AB_{0}B_{1}(w)\,F_{D}^{(r-1)}\left( \frac{\nu +|m|}{2};\frac{m_{1}}{2}%
,\cdots ,\frac{m_{r-1}}{2};\frac{|m|}{2};t(w)u_{1},\cdots
,t(w)u_{r-1}\right) ,  \label{FD}
\end{gather}
where $F_{D}$ is a Lauricella function (Srivastava and Karlsson (1985, p.
41) where we correct the typographical error with Equation \ref{FD})).
Equation \ref{FD} gives a representation of the $pdf$ of the distribution $W.
$ We will show in Theorem \ref{th2} that the $cdf$ of $W$ is also a
Lauricella $F_{D}^{(r)}$ function. This representation will yield a
numerically computable algorithm for finding $p$-values. Equation \ref{hw
line 2} yields a numerically tractable expression for the $pdf$ of $W$. In
Section \ref{local bound pdf}, we give a tight local truncation error bound $%
e_{\tau }(w)$ for determining the number of terms $\tau $ to use in the
partial sum expression.

\subsection{Exact Expressions for the $pdf$ of $W$ with $r=2$}

If $r=2$, Equation \ref{Wline3} is a hypergeometric series, and we have the
following result.

\begin{theorem}
With the notation above, $a=\nu \alpha _{2}/|m|,$ and $r=2,$ the $pdf$ of $W$
is given by 
\begin{eqnarray}
h_{W}(w) &=&AB_{0}B_{1}(w)\sum_{j=0}^{\infty }\frac{\left( \frac{\nu
+m_{1}+m_{2}}{2}\right) _{j}\left( \frac{m_{1}}{2}\right) _{j}}{j!\left( 
\frac{m_{1}+m_{2}}{2}\right) _{j}}(u_{1}t(w))^{j}  \label{line1ofbigth} \\
&=&AB_{0}B_{1}(w)\,_{2}F_{1}\left( \QDATOP{\frac{\nu +m_{1}+m_{2}}{2},\frac{%
m_{1}}{2}}{\frac{m_{1}+m_{2}}{2}};(1-\frac{\alpha _{r}}{\alpha _{1}})\frac{w%
}{a+w}\right) .  \label{bigth}
\end{eqnarray}
\end{theorem}

To find the $cdf$ of $W$ when $r=2$, integrate $h_{W}(w)$ in Equation \ref
{bigth}.

We note that if we had used the notation of Kotz, Johnson, and Boyd (1967)
and scaled $y$ by $y/\delta $ with $0<\delta <\alpha _{r}$, then $u_{2}>0$.
In this situation, we would use the Bailey transformation (Srivastava and
Karlsson (1985, p. 304)) to convert the two variable hypergeometric series
in Equation \ref{line1ofbigth} to the $_{2}F_{1}$ function in Equation \ref
{bigth}.

\subsection{Exact Expressions for the $cdf$ of $W$ with $r\geq 2$}

The Lauricella function $F_{D}^{(r-1)}$ in Equation \ref{FD} has an integral
representation (Exton, 1976, p. 49) where the domain of integration is over
the simplex $E_{r}$ with $x_{1}+\cdots +x_{r}=1\,(x_{i}\geq 0,1\leq i\leq r)$
as 
\begin{eqnarray}
&&F_{D}^{(r-1)}\left( \frac{\nu +|m|}{2};\frac{m_{1}}{2},\cdots ,\frac{%
m_{r-1}}{2};\frac{|m|}{2};t(w)u_{1},\cdots ,t(w)u_{r-1}\right)   \notag \\
&=&\Gamma \left[ \frac{\frac{|m|}{2}}{\frac{m_{1}}{2},\cdots ,\frac{m_{r}}{2}%
}\right] \int_{E_{r}}(1-\sum_{i=1}^{r-1}t(w)u_{i}x_{i})^{-\frac{\nu +|m|}{2}%
}\prod_{i=1}^{r}x_{i}^{\frac{m_{i}}{2}-1}d\mathbf{x.}  \label{muti-integral}
\end{eqnarray}
In Dunkl and Ramirez (1994a, 1994b), we computed the surface measure of
ellipsoids using hyperelliptic integrals. We showed that the $(n-1)$%
-dimensional hyperelliptic integral could be transformed into a univariate
integral using the Euler integral representation (Exton, 1976, p. 49) for $%
F_{D}.$ This transformation does not apply to Equation \ref{muti-integral}
since $\frac{\nu +|m|}{2}>\frac{|m|}{2}.$ Here we will use a different
approach.

We show how to represent the $cdf$ of the generalized $F$ distribution $W$
as a Lauricella $F_{D}^{(r)}$ function. This representation will provide a
numerically tractable procedure for computing the $cdf$ of $W$, denoted by $%
H_{W}(w),$ which does not require integrating the $pdf$ of $W.$

\begin{theorem}
\label{th2}With the notation above and $r\geq 2,$ the $cdf$ of $W$ is given
by 
\begin{gather}
H_{W}(y)=AB_{0}\frac{y^{|m|/2}}{(|m|/2)(a+y)^{(\nu +|m|)/2}}  \notag \\
F_{D}^{(r)}(\frac{\nu +|m|}{2};\frac{m_{1}}{2},\cdots ,\frac{m_{r-1}}{2},1;%
\frac{|m|}{2}+1;t(y)u_{1},\cdots ,t(y)u_{r-1},t(y)),  \label{HW}
\end{gather}
with $a=\nu \alpha _{r}/|m|$ and $t(y)=y/(a+y)$ as before.
\end{theorem}

\begin{proof}
From Equations \ref{FD} and \ref{muti-integral}, write the $cdf$ of $W$ as 
\begin{gather*}
H_{W}(y)=\int_{0}^{y}h_{W}(w)dw \\
=AB_{0}\frac{\Gamma \left( \frac{|m|}{2}\right) }{\prod_{i=1}^{r}\Gamma
\left( \frac{m_{i}}{2}\right) }\int_{0}^{y}\frac{w^{|m|/2-1}}{(a+w)^{(\nu
+|m|)/2}} \\
\int_{E_{r}}\left( 1-\sum_{i=1}^{r}\frac{w}{a+w}u_{i}x_{i}\right) ^{-(\nu
+|m|)/2}\prod_{i=1}^{r}x_{i}^{m_{i}/2-1}d\mathbf{x}dw \\
=AB_{0}\frac{\Gamma \left( \frac{|m|}{2}\right) }{\prod_{i=1}^{r}\Gamma
\left( \frac{m_{i}}{2}\right) }\int_{0}^{y}\int_{E_{r}}w^{r-1}(a+w-%
\sum_{i=1}^{r}wu_{i}x_{i})^{-(\nu +|m|)/2} \\
\prod_{i=1}^{r}(wx_{i})^{m_{i}/2-1}d\mathbf{x}dw.
\end{gather*}
Change variables with $s_{i}=wx_{i}/y$ $(1\leq i\leq r)$ and $s_{r+1}=1-w/y.$
Note that $\sum_{i=1}^{r}s_{i}=w/y$ with the absolute value of the inverse
Jacobian $J^{-1}=\left| \frac{\partial (x_{1},\cdots ,x_{r-1},w)}{\partial
(s_{1},\cdots ,s_{r-1},s_{r})}\right| =w^{r-1}/y^{r}.$ Thus 
\begin{gather}
H_{W}(y)=AB_{0}\frac{\Gamma \left( \frac{|m|}{2}\right) }{%
\prod_{i=1}^{r}\Gamma \left( \frac{m_{i}}{2}\right) }  \notag \\
\int_{E_{r+1}}y^{r}\left( a+y(1-s_{r+1})-\sum_{i=1}^{r}ys_{i}u_{i}\right)
^{-(\nu +|m|)/2}\prod_{i=1}^{r}(ys_{i})^{m_{i}/2-1}d\mathbf{s}  \notag \\
=AB_{0}\frac{\Gamma \left( \frac{|m|}{2}\right) }{\prod_{i=1}^{r}\Gamma
\left( \frac{m_{i}}{2}\right) }(a+y)^{-(\nu +|m|)/2}y^{|m|/2}dw  \notag \\
\int_{E_{r+1}}\left( 1-\frac{ys_{r+1}}{a+y}-\sum_{i=1}^{r}\frac{y}{a+y}%
s_{i}u_{i}\right) ^{-(\nu +|m|)/2}\prod_{i=1}^{r}s_{i}{}^{m_{i}/2-1}d\mathbf{%
s}  \notag \\
=AB_{0}\frac{1}{|m|/2}\frac{y^{|m|/2}}{(a+y)^{(\nu +|m|)/2}}\,
\label{FD_r>2} \\
F_{D}^{(r)}(\frac{\nu +|m|}{2};\frac{m_{1}}{2},\cdots ,\frac{m_{r-1}}{2},1;%
\frac{|m|}{2}+1;t(y)u_{1},\cdots ,t(y)u_{r-1},t(y)),  \notag
\end{gather}
with $a=\nu \alpha _{r}/|m|$ and $t(y)=y/(a+y).$
\end{proof}

To convert Equation \ref{FD_r>2} into a numerically tractable series, write 
\begin{gather}
H_{W}(y)=AB_{0}\frac{1}{|m|/2}\frac{y^{|m|/2}}{(a+y)^{(\nu +|m|)/2}}\, 
\notag \\
F_{D}^{(r)}(\frac{\nu +|m|}{2};\frac{m_{1}}{2},\cdots ,\frac{m_{r-1}}{2},1;%
\frac{|m|}{2}+1;t(y)u_{1},\cdots ,t(y)u_{r-1},t(y))  \notag \\
=B_{0}\frac{y^{|m|/2}}{(a+y)^{(\nu +|m|)/2}}  \notag \\
\sum_{j=0}^{\infty }\frac{(\frac{\nu +|m|)}{2})_{j}}{(\frac{|m|}{2})_{j+1}}%
\left( \frac{y}{a+y}\right) ^{j}\left[ A\sum_{i_{1}+\cdots +i_{r-1}\leq
j}\prod_{k=1}^{r-1}\frac{u_{k}^{i_{k}}}{i_{k}!}\left( \frac{m_{k}}{2}\right)
_{i_{k}}\right]  \notag \\
=B_{0}\frac{y^{|m|/2}}{(a+y)^{(\nu +|m|)/2}}\sum_{j=0}^{\infty }\frac{(\frac{%
\nu +|m|)}{2})_{j}}{(\frac{|m|}{2})_{j+1}}\left( \frac{y}{a+y}\right)
^{j}\sum_{i=0}^{j}c_{i}.  \label{useful}
\end{gather}

\section{Local Truncation Error Bounds}

Denote by $\widehat{h}_{W}(w)$ and $\widehat{H}_{W}(y)$ the partial sum
estimates for $h_{W}(w)$ and$\ H_{W}(y)$, respectively, from Equations \ref
{hw line 2} and \ref{HW line 2}. In this Section, we derive local truncation
error bounds to determine the number of terms required by the partial sums.

\subsection{Local Truncation Error Bound $e_{\protect\tau }^{\ast }(y)$ for
the $cdf$ of $W$\label{bound section}}

For Equation \ref{useful} to be numerically tractable, we derive the local
truncation error. Write $t(y)=(y/(a+y))<1,$ 
\begin{gather}
H_{W}(y)=B_{0}yB_{1}(y)\frac{y}{(|m|/2)}\sum_{j=0}^{\infty }\frac{(\frac{\nu
+|m|)}{2})_{j}}{(\frac{|m|}{2}+1)_{j}}t(y)^{j}\sum_{i=0}^{j}c_{i}  \notag \\
=B_{0}B_{1}(y)\frac{y}{(|m|/2)}\sum_{j=0}^{\tau }\frac{(\frac{\nu +|m|)}{2}%
)_{j}}{(\frac{|m|}{2}+1)_{j}}t(y)^{j}\sum_{i=0}^{j}c_{i}+  \label{HW line 2}
\\
B_{0}B_{1}(y)\frac{y}{(|m|/2)}\sum_{j=\tau +1}^{\infty }\frac{(\frac{\nu
+|m|)}{2})_{j}}{(\frac{|m|}{2}+1)_{j}}t(y)^{j}\sum_{i=0}^{j}c_{i}  \notag
\end{gather}

\begin{gather}
=\widehat{H}_{W}(y)+B_{0}B_{1}(y)\frac{y}{(|m|/2)}\frac{(\frac{\nu +|m|)}{2}%
)_{\tau +1}}{(\frac{|m|}{2}+1)_{\tau +1}}t(y)^{\tau +1}\times  \notag \\
\sum_{j=0}^{\infty }\frac{(\frac{\nu +|m|)}{2}+\tau +1)_{j}}{(\frac{|m|}{2}%
+\tau +2)_{j}}t(y)^{j}(1-1+\sum_{i=0}^{\tau +1+j}c_{i}) \\
=\widehat{H}_{W}(y)+B_{0}B_{1}(y)\frac{y}{(|m|/2)}\frac{(\frac{\nu +|m|)}{2}%
)_{\tau +1}}{(\frac{|m|}{2}+1)_{\tau +1}}t(y)^{\tau +1}\,_{2}F_{1}\left( 
\QATOP{\frac{\nu +|m|)}{2}+\tau +1,1}{\frac{|m|}{2}+\tau +2};t(y)\right) 
\notag \\
-B_{0}B_{1}(y)\frac{y}{(|m|/2)}\frac{(\frac{\nu +|m|)}{2})_{\tau +1}}{(\frac{%
|m|}{2}+1)_{\tau +1}}t(y)^{\tau +1}\sum_{j=0}^{\infty }\frac{(\frac{\nu +|m|)%
}{2})_{j}}{(\frac{|m|}{2}+1)_{j}}t(y)^{j}(1-\sum_{i=0}^{\tau +1+j}c_{i}). 
\notag
\end{gather}

The partial sum estimate $\widehat{H}_{W}(y)$ can be enhanced by identifying
most of the truncation error as a scaled $_{2}F_{1}$ hypergeometric
functions. The remaining truncation error is bounded by a scaled $_{2}F_{1}$
function and is stated in the following.

\begin{theorem}
With the notation above, the estimated $P[W\leq y]$ is given by 
\begin{equation}
\widehat{H}_{W}(y)+B_{0}B_{1}(y)\frac{y}{(|m|/2)}\frac{(\frac{\nu +|m|)}{2}%
)_{\tau +1}}{(\frac{|m|}{2}+1)_{\tau +1}}t(y)^{\tau +1}\,_{2}F_{1}\left( 
\QATOP{\frac{\nu +|m|)}{2}+\tau +1,1}{\frac{|m|}{2}+\tau +2};t(y)\right)
\label{H estimate}
\end{equation}
with local truncation error bound given by 
\begin{gather}
e_{\tau }^{\ast }(y)=(1-\sum_{i=0}^{\tau +1}c_{i})B_{0}B_{1}(y)\frac{y}{%
(|m|/2)}\frac{(\frac{\nu +|m|)}{2})_{\tau +1}}{(\frac{|m|}{2}+1)_{\tau +1}}%
t(y)^{\tau +1}\times \,  \label{H error2} \\
_{2}F_{1}\left( \QATOP{\frac{\nu +|m|)}{2}+\tau +1,1}{\frac{|m|}{2}+\tau +2}%
;t(y)\right) .  \notag
\end{gather}
\end{theorem}

To find $\tau $, we increase the size of $\tau $ unless the remaining error $%
e_{\tau }^{\ast }(y)$ from Equation \ref{H error2} is less than a prescribed
small value. The suggested value is $10^{-4}.$

\subsection{Local Truncation Error Bound $e_{\protect\tau }(w)$ for the $pdf$
of $W\label{local bound pdf}$}

Recall that Equation \ref{hw line 2} yields a numerically tractable
expression for the $pdf$ of $W$. A tight local truncation error bound for
determining the number of terms $\tau $ to use in the partial sum expression
follows as above and is stated in the following.

\begin{theorem}
With the notation above, 
\begin{equation}
\widehat{h}_{W}(w)=B_{0}B_{1}(w)\sum_{j=0}^{\tau }\frac{c_{j}\left( \frac{%
\nu +|m|}{2}\right) _{j}}{\left( \frac{|m|}{2}\right) _{j}}t(w)^{j}
\label{hW error}
\end{equation}
with local truncation error bound given by 
\begin{equation}
e_{\tau }(w)=c_{\tau +1}B_{0}B_{1}(w)\frac{(\frac{\nu +|m|)}{2})_{\tau +1}}{(%
\frac{|m|}{2})_{\tau +1}}t(w)^{\tau +1}\,_{2}F_{1}\left( \QATOP{\frac{\nu
+|m|)}{2}+\tau +1,1}{\frac{|m|}{2}+\tau +2};t(w)\right) ,
\label{local trun error for h}
\end{equation}
a scaled $_{2}F_{1}$ hypergeometric function.
\end{theorem}

To determine the number of terms for the partial sum estimate $\widehat{h}%
_{W}(w),$ increase the size of $\tau $ unless the local truncation error
from Equation \ref{local trun error for h} is less than a prescribed small
value. The suggested value is $y10^{-4}$ where the $p$-value is calculated
from $y.$

\section{Applications}

We will give two applications where the distribution of the test statistic
is the generalized $F$ distribution.

\subsection{Detection of Outliers}

Cook's (1977) $D_{I}$ statistics are used widely for assessing influence of
design points in regression diagnostics. These statistics typically contain
a leverage component and a standardized residual component. Subsets having
large $D_{I}$ are said to be influential, reflecting high leverage for these
points or that $I$ contains some outliers from the data. Consider the linear
model 
\begin{equation}
\mathbf{Y}_{0}=\mathbf{X}_{0}\mathbf{\beta }+\mathbf{\varepsilon }_{0},
\label{linmodel}
\end{equation}
where $\mathbf{Y}_{0}$ is a $(N\times 1)$ vector of observations, $\mathbf{X}%
_{0}$ is a $(N\times k)$ full rank matrix of known constants, $\mathbf{\beta 
}$ is a $(k\times 1)$ vector of unknown parameters, and $\mathbf{\varepsilon 
}_{0}$ is a $(N\times 1)$ vector of randomly distributed Gaussian errors
with $E(\mathbf{\varepsilon }_{0})=\mathbf{0}$ and $Var(\mathbf{\varepsilon }%
_{0})=\sigma ^{2}\mathbf{I}_{N}.$ The least squares estimate of $\mathbf{%
\beta }$ is $\mathbf{\hat{\beta}}=(\mathbf{X}_{0}^{^{\prime }}\mathbf{X}%
_{0})^{-1}\mathbf{X}_{0}^{^{\prime }}\mathbf{Y}_{0}.$ The basic idea in 
\textit{influence analysis,} as introduced by Cook (1977), concerns the
stability of a linear regression model under small perturbations. For
example, if some cases are deleted, then what changes occur in estimates for
the parameter vector $\mathbf{\beta }?$ Cook's $D_{I}$ statistics are based
on a Mahalanobis distance between $\mathbf{\hat{\beta}}$ (using all the
cases) and $\mathbf{\hat{\beta}}_{I}$ (using all cases except those in the
subset $I)$, as given by 
\begin{equation}
D_{I}(\mathbf{\hat{\beta}},\mathbf{M},c\hat{\sigma}^{2})=(\mathbf{\hat{\beta}%
}_{I}-\mathbf{\hat{\beta}})^{^{\prime }}\mathbf{M}(\mathbf{\hat{\beta}}_{I}-%
\mathbf{\hat{\beta}})/(c\hat{\sigma}^{2}),  \label{cookd}
\end{equation}
with a $(k\times k)$ nonnegative definite matrix $M$, $\hat{\sigma}^{2}$ is
an unbiased estimate of the variance, and a user defined constant $c$. We
use $c=r$ and the estimator $s_{I}^{2}$, the sample variance estimator with
the cases in $I$ omitted We will discuss the case with $\mathbf{M}=$ $%
\mathbf{X}^{^{\prime }}\mathbf{X}$, where $\mathbf{X}$ denotes the remaining
rows of $\mathbf{X}_{0}$. We have chosen $s_{I}^{2}$ as the estimator for $%
\sigma ^{2}$ since this estimator and the numerator of Equation \ref{cookd}
are independent.

Using the results in this paper, we are able to numerically compute the $cdf$
of Cook's $D_{I}$ statistics in the case of joint outliers, and, in
particular, to compute the $p$-values for $D_{I}$. This approach provides a
statistical procedure for identifying influential observations based on $p$%
-values.

\subsubsection{\textbf{Notation}}

To fix the notation, let $I$ be a subset of $\{1,\ldots ,N\},$ say $%
I=\{i_{1},\ldots ,i_{r}\}.$ Let $\mathbf{X}_{0}$ be partitioned as $\mathbf{X%
}_{0}^{^{\prime }}=[\mathbf{X}^{^{\prime }},\mathbf{Z}^{^{\prime }}],$ with $%
\mathbf{X}$ containing the rows determined by $I$, and $\mathbf{Z}$ the
remaining rows. We assume that the matrices $\mathbf{X}_{0}$, $\mathbf{X}$,
and $\mathbf{Z}$ all of full rank, of orders $(N\times k)$, $(n\times k)$,
and $(r\times k)$, respectively such that $k<n<N$, and $n+r=N,$ with $r<k$
for notational convenience. Partition $\mathbf{Y}_{0}^{\prime }=[\mathbf{Y}%
_{1}^{\prime },\mathbf{Y}_{2}^{\prime }]$, and $\mathbf{\varepsilon }%
_{0}^{^{\prime }}=[\mathbf{\varepsilon }_{1}^{^{\prime }},\mathbf{%
\varepsilon }_{2}^{^{\prime }}]$. Thus Equation \ref{linmodel} has been
transformed into 
\begin{equation}
\left[ 
\begin{tabular}{r}
$\mathbf{Y}_{1}$ \\ 
$\mathbf{Y}_{2}$%
\end{tabular}
\right] =\left[ 
\begin{tabular}{r}
$\mathbf{X}$ \\ 
$\mathbf{Z}$%
\end{tabular}
\right] \mathbf{\beta }+\left[ 
\begin{tabular}{r}
$\mathbf{\varepsilon }_{1}$ \\ 
$\mathbf{\varepsilon }_{2}$%
\end{tabular}
\right] .  \label{block_model}
\end{equation}

The ordered eigenvalues of $\mathbf{Z}(\mathbf{X}_{0}^{^{\prime }}\mathbf{X}%
_{0})^{-1}\mathbf{Z}^{^{\prime }}$ are denoted $\{\lambda _{1}\geq \cdots
\geq \lambda _{r}>0\}$ usually called the canonical leverages. Jensen and
Ramirez (1991) showed that the $cdf$ for $W_{0}=T/V,$ equivalently for $%
W=(T/r)/(V/\nu ),$ is a weighted series of $F$ distributions, and they
computed the stochastic bounds 
\begin{equation}
F_{r}(w\,;\alpha _{1};\nu )\leq F_{r}(w\,;\alpha _{1},\ldots ,\alpha
_{r};1,\ldots ,1;\nu )\leq F_{r}(w\,;\alpha ^{\ast };\nu )\ ,  \label{bounds}
\end{equation}
with the maximum weight $\alpha _{1}$, the geometric mean $\alpha ^{\ast }$
of the weights $\{\alpha _{1,}\ldots ,\allowbreak \alpha _{r}\},$ and $%
F_{r}(w\,;\alpha ;\nu )$ the scaled central $F$ distribution.

The basic characterization theorem for $D_{I}$ is given in Jensen and
Ramirez (1998a) and is:

\begin{theorem}
\label{basictheorem}Suppose that $\mathcal{L}(\mathbf{Y})=N_{N}(\mathbf{X}%
_{0}\mathbf{\beta },\sigma ^{2}\mathbf{I}_{N})$, then the distribution of $%
D_{I}(\mathbf{\hat{\beta}},\mathbf{X}^{^{\prime }}\mathbf{X},rs_{I}^{2})$ is
given by $F_{r}(w;\allowbreak \lambda _{1},\cdots ,\lambda _{r};1,\cdots
,1;N-r-k)$.
\end{theorem}

With $r=1$, $\mathcal{L}(D_{i}(\mathbf{\hat{\beta}},\mathbf{X}^{^{\prime }}%
\mathbf{X},s_{i}^{2})/\lambda _{1})=F(1,N-1-k)$. Outliers also can be tested
using the studentized deleted residuals with $\mathcal{L}((y_{i}-\hat{y}%
_{(i)})/(s_{i}(1+\allowbreak \mathbf{x}_{i}(\mathbf{X}^{^{\prime }}\mathbf{X}%
)^{-1}\mathbf{x}_{i}^{^{\prime }})^{1/2}))=t(N-1-k)$ where $\hat{y}_{(i)}$
denotes the predicted value using $(\mathbf{Y}_{1},\mathbf{X});$ or with the
externally studentized residuals ($RStudent$) with $\mathcal{L(}(y_{i}-\hat{y%
}_{i})/\allowbreak (s_{i}\sqrt{1-h_{ii}}))=t(N-1-k)$ where $\hat{y}_{i} $
denotes the predicted value using $(\mathbf{Y},\mathbf{X}_{0})$ and $h_{ii} $
is the canonical leverage also denoted as $\lambda _{1}$. In Jensen and
Ramirez (1998b) it is shown that the $p$-values from these two tests are
also equal to the $p$-values from Theorem \ref{basictheorem}. Thus, in case
of single deletion with $r=1$, all of these three standard tests for
outliers will have a common $p$-value.

\subsubsection{Examples}

For the Hald (1952, p. 647) data set ($N=13$ and $k=5)$ using the test
statistic $D_{I}(\mathbf{\hat{\beta}},\mathbf{X}^{^{\prime }}\mathbf{X}%
,2s_{I}^{2})$ and the global bounds in Equation \ref{bounds}, we can show
that the only pair $(r=2)$ of observations (from the 78 possible pairs)
which could possibly be influential at the 5\% significance level is $%
I=\{6,8\}$ with $0.01305<p_{I}<0.04610$. Using $m_{1}=m_{2}=1,$ the
canonical leverages $\mathbf{\lambda }=(0.408676$, $0.124019)$ for the
weights $\mathbf{\alpha }$, the degrees of freedom $\nu =N-r-k=6$, and the
observed Cook's $D_{I}$ statistic $y=2.19331$, we can now easily compute
from Equation \ref{bigth} that the $p$-value is $p_{I}$= $0.02181$.

For the Longley (1967) data set, Cook (1977) noted that observations 5 and
16 may be influential. To test for the joint influence of $I=\{5,16\}$, we
use the test statistic $D_{I}(\mathbf{\hat{\beta}},\mathbf{X}^{^{\prime }}%
\mathbf{X},2s_{I}^{2}),$ with $r=2$, the canonical leverages $\mathbf{%
\lambda }=(0.690029,0.614130)$ for the weights, $\nu =N-r-k=16-2-7=7$, and
the observed Cook's $D_{I}$ statistic $y=1.812433$, we compute that the $p$%
-value is $p_{I}=$ $0.12927$.

Using the test statistic $D_{I}(\mathbf{\hat{\beta}},\mathbf{X}^{^{\prime }}%
\mathbf{X},2s_{I}^{2})$ and the global bounds Equation \ref{bounds}, it is
easy to compute that the only possible pairs that need to be considered at
the $5\%$ significance level are (1) $I_{1}=\{4,5\}$ with $\mathbf{\lambda }%
=(0.615959,0.371827)$, $y=2.57861$, and $0.03822\leq p_{I_{1}}=0.04186\leq
0.06356,$ (2) $I_{2}=\{4,15\}$ with $\mathbf{\lambda }=(0.505387.0.393672)$, 
$y=1.76885$, and $0.04961\leq p_{I_{2}}=0.04982\leq 0.05555$, and (3) $%
I_{3}=\{10,16\}$ with $\mathbf{\lambda }=(0.736874,0.695572)$, $y=2.57906$,
and $0.03761\leq p_{I_{3}}=0.04571\leq 0.07979$ where the $p$-values $p_{I}$
are computed from Equation \ref{bigth}$.$

Our recommendation to the practitioner, who wishes to find joint outliers,
is to initially screen for potential joint outliers using Equation \ref
{bounds} with $D_{I}(\mathbf{\hat{\beta}},\mathbf{X}^{^{\prime }}\mathbf{X}%
,rs_{I}^{2})$. If $r=1$ then the distribution of $D_{i}$ is a scaled central 
$F$ distribution. If $r=2$ then the distribution of $D_{I}$ is a scaled $%
_{2}F_{1}$ series. If $r>2$ then use Equation\ref{H error2} to find the
numbers of terms required to have the local truncation error small. The
suggested value for the bound is $10^{-4}.$ The $p$-values for the $cdf$ for
the distribution of $D_{I}(\mathbf{\hat{\beta}},\mathbf{X}^{^{\prime }}%
\mathbf{X},\allowbreak rs_{I}^{2})$ are calculated using the enhanced
truncated series in Equation \ref{H estimate}.

\subsection{Misspecified Hotelling's $T$ test}

Hotelling's $T^{2}$ is used widely in multivariate data analysis,
encompassing tests for means, the construction of confidence ellipsoids, the
analysis of repeated measurements, and statistical process control. To
support a knowledgeable use of $T^{2}$, its properties must be understood
when model assumptions fail. Jensen and Ramirez (1991) have studied the
misspecification of location and scale in the model for a multivariate
experiment under practical circumstances to be described.

To set the notation, let $N_{p}(\mathbf{\mu },\mathbf{\Sigma })$ be the
Gaussian distribution with mean $\mathbf{\mu }$, and dispersion $\mathbf{%
\Sigma }$ and let $W_{p}(\nu ^{\ast },\mathbf{\Sigma })$ denote the central
Wishart distribution having $\nu ^{\ast }$ degrees of freedom and scale
parameter $\mathbf{\Sigma }$.\ Consider the representation $T^{2}=\nu ^{\ast
}\mathbf{Y}^{\prime }\mathbf{W}^{-1}\mathbf{Y}$ where $(\mathbf{Y},\mathbf{W}%
)$ are independent and $\mathcal{L}(\mathbf{Y})=N_{p}(\mathbf{\mu },\mathbf{%
\Sigma })$ as before, but now $\mathcal{L}(\mathbf{W})=W_{p}(\mathbf{\nu }%
^{\ast },\mathbf{\Omega })$. Denote the ordered roots of $\mathbf{\Omega }^{-%
\frac{1}{2}}\mathbf{\Sigma \Omega }^{-\frac{1}{2}}$ by $\{\pi _{1}\geq \pi
_{2}\geq \cdots \geq \pi _{p}>0\}$. A principal result for $T^{2}$ under
misspecified scale is given in Jensen and Ramirez (1991) and is the
following.

\begin{theorem}
The distribution of the test statistic $((\nu ^{\ast }-p+1)/p)(T^{2}/\nu
^{\ast })$ is the generalized $F$ distribution $F_{r}(w;\pi _{1},\allowbreak
\cdots ,\allowbreak \pi _{p};1,\cdots ,1;\nu ^{\ast }-p+1)$.
\end{theorem}

\subsubsection{Hotelling's misspecifed scale distribution}

The conventional model for $T^{2}$ is based on a random sample $%
\{X_{1},\ldots ,X_{N}\}$ from $N_{p}(\mathbf{\mu },\mathbf{\Sigma })$ using
the unbiased sample means and dispersion matrix $(\mathbf{\bar{X}},\mathbf{S}%
).$ We have $\mathcal{L}(\mathbf{\bar{X})=N}_{p}(\mathbf{\mu },\frac{1}{N}%
\mathbf{\Sigma })$ and $\mathcal{L}((N-1)\mathbf{S)=}W_{p}(N-1,\mathbf{%
\Sigma })$, or $\mathcal{L}(\frac{N-1}{N}\mathbf{S)=}W_{p}(N-1,\frac{1}{N}%
\mathbf{\Sigma }).$ Thus $T^{2}=(N-1)(\mathbf{\bar{X}-\mu })^{\prime }(\frac{%
N-1}{N}\mathbf{S)}^{-1}(\mathbf{\bar{X}-\mu })=N_{p}(\mathbf{\bar{X}-\mu }%
)^{\prime }\mathbf{S}^{-1}(\mathbf{\bar{X}-\mu })$ and $\mathcal{L}%
(((N-p)/p)(T^{2}/(N-1)))=F(p,N-p)$, the central $F$ distribution when $N>p$.
If the process dispersion parameters have shifted, then $T^{2}$ is
misspecified with $\mathcal{L}((N-1)\mathbf{S)=}W_{p}(N-1,\mathbf{\Omega })$%
, and with $((N-p)/p)(T^{2}/(N-1))$ the generalized $F$ distribution $%
F_{r}(w;\pi _{1},\cdots ,\pi _{p};1,\cdots ,1;N-p)$. Here $r=p,$ $\nu =\nu
^{\ast }-p+1=N-p$, and $\{\pi _{1}\geq \pi _{2}\geq \cdots \geq \pi _{p}>0\}$
the ordered roots of $\mathbf{\Omega }^{-\frac{1}{2}}\mathbf{\Sigma \Omega }%
^{-\frac{1}{2}}$.

\subsubsection{Examples}

An important application of generalized $F$ distributions is for computing
the power of a misspecified Hotelling's $T^{2}$ test for a multivariate
quality control chart. Power analysis for a misspecified mean $\mathbf{\mu }$
is standard. Using generalized $F$ distributions, the power analysis for a
misspecified covariance $\mathbf{\Omega }$ can be performed. If a process
changes, not only will the mean change but generally the covariance
structure will also change. The robustness of $T^{2}$ under misspecification
of scale can be verified by computing the cumulative density of $T^{2}$ for
varying choices of $\pi _{1}\geq \pi _{2}\geq \cdots \geq \pi _{p}>0$ at the
critical value of $T^{2}$. For example, if $\mathbf{\Omega }_{\rho }$ is a $%
3\times 3$ equicorrelated matrix ($r=p=3$) with $\rho =0.5$, and if $\mathbf{%
\Sigma }$ is the identity matrix, then the eigenvalues of\textbf{\ }$\Omega
_{\rho }^{-\frac{1}{2}}\mathbf{\Sigma \Omega }_{\rho }^{-\frac{1}{2}}$ are $%
\{\pi _{1}=(1-\rho )^{-1},\pi _{2}=(1-\rho )^{-1},\pi _{3}=(1+2\rho
)^{-1}\}=\{2,2,1/2\}$. If $N=12$ with $\nu =N-p=9$, the nominal 95\%
critical value of $((N-p)/p)(T^{2}/(N-1))$ is $F(0.95;p,N-p)=3.8625$.
However, the exact right-hand tail probability for $Y=((N-p)/p)(T^{2}/(N-1))$
is not $0.05$ but rather $P[Y=((N-p)/p)(T^{2}/(N-1))\geq 3.8625]=0.12310$.
In this example, $\pi _{1}=\pi _{2},$ so we could compute the $p$-values
exactly from Theorem 1, with $F_{3}(w;\pi _{1},\pi _{2},\pi
_{3};1,1,1;N-p)=F_{2}(w;\pi _{1},\pi _{3};2,1;N-p).$ Instead, we use this
problem to demonstrate the number of terms required by the three numerical
methods discussed in this paper.

In Table 1, we present similar computations for varying $\rho $. For each $%
\rho $ in the Table 1, and with the corresponding eigenvalues $\pi _{1}\geq
\pi _{2}\geq \pi _{3}>0$ of $\mathbf{\Omega }_{\rho }^{-\frac{1}{2}}\mathbf{%
\Sigma \Omega }_{\rho }^{-\frac{1}{2}},$ we give the value of $%
P[Y=((N-p)/p)(T^{2}/(N-1))\geq 3.8625.$ Also shown are the number of terms
required using the three numerical presented in this paper. The first is $%
\tau _{1}$ from Equation \ref{gWerr} required to satisfy $ye_{\tau _{1}}\leq
10^{-4},$ the second is $\tau _{2}$ from Equation \ref{local trun error for
h} required to satisfy $ye_{\tau _{2}}(y)\leq 10^{-4},$ and the third is $%
\tau _{3}$ from Equation \ref{H error2} required to satisfy $e_{\tau
_{2}}^{\ast }(y)\leq 10^{-4}$. The inputs are $r=3$, the weights $\pi
_{1}\geq \pi _{2}\geq \pi _{3}>0$, $\nu =N-p=12-3=9,$ and $y=3.8625$.

\begin{center}
\begin{tabular}{|c|c|c|c|c|}
\multicolumn{5}{c}{\textbf{Table 1. Misspecified Type I Error}} \\%
[1mm] \hline
$\rho $ & $\tau _{1}$ & $\tau _{2}$ & $\tau _{3}$ & $P[Y\geq 3.8625]$ \\ 
\hline
0.0 & 1 & 1 & 1 & 0.0500 \\ 
0.1 & 6 & 7 & 6 & 0.0526 \\ 
0.2 & 10 & 11 & 8 & 0.0600 \\ 
0.3 & 15 & 15 & 12 & 0.0727 \\ 
0.4 & 20 & 20 & 16 & 0.0926 \\ 
0.5 & 28 & 26 & 21 & 0.1231 \\ 
0.6 & 40 & 32 & 27 & 0.1704 \\ 
0.7 & 58 & 40 & 34 & 0.2458 \\ 
0.8 & 92 & 49 & 43 & 0.3712 \\ 
0.9 & 185 & 58 & 55 & 0.59055 \\ \hline
\end{tabular}
\end{center}

As anticipated, the numbers of terms $\tau $ required is fewer when the
enhanced partial sum from Equation \ref{H estimate} is used. More
importantly, the method from Section \ref{bound section} does not require
that the $pdf$ to be numerically integrated.

\section{Conclusion}

We have derived the exact distribution of the generalized $F$ distribution%
\newline
$F_{2}(w;\allowbreak \alpha _{1},\allowbreak \alpha _{2};\allowbreak
m_{1},m_{2};2)$ in terms of the hypergeometric series $_{2}F_{1}.$ This
extends the corresponding result of Bock and Solomon for a mixture of two
chi-square distributions to the generalized $F$ distribution with $r=2.$
Explicit representations for the case $r\geq 2$ are given in terms of a
Lauricella $F_{D}$ functions. Numerically computable series expansion have
been derived. Applications to the detection of joint outliers and to the
misspecified Hotelling $T^{2}$ statistic have been given.

\end{document}